\documentclass[11pt]{article}
\usepackage{amsmath, upgreek}
\usepackage{amssymb}
\usepackage{color}
\usepackage{amscd}
\usepackage{xspace}
\usepackage{verbatim}
\usepackage{graphicx}
\renewcommand{\theequation}{\thesection.\arabic{equation}
}
\newcommand{\mysection}[1]{
\section{#1}\setcounter{equation}{0}}
 \usepackage{cite}
 
\begin{document}


\newcommand{\txt}[1]{\;\text{ #1 }\;}
\newcommand{\tbf}{\textbf}
\newcommand{\tit}{\textit}
\newcommand{\tsc}{\textsc}
\newcommand{\trm}{\textrm}
\newcommand{\mbf}{\mathbf}
\newcommand{\mrm}{\mathrm}
\newcommand{\bsym}{\boldsymbol}
\newcommand{\scs}{\scriptstyle}
\newcommand{\sss}{\scriptscriptstyle}
\newcommand{\txts}{\textstyle}
\newcommand{\dsps}{\displaystyle}
\newcommand{\fnz}{\footnotesize}
\newcommand{\scz}{\scriptsize}
\newcommand{\be}{\begin{equation}}
\newcommand{\bel}[1]{\begin{equation}\label{#1}}
\newcommand{\ee}{\end{equation}}
\newcommand{\eqnl}[2]{\begin{equation}\label{#1}{#2}\end{equation}}
\newcommand{\barr}{\begin{eqnarray}}
\newcommand{\earr}{\end{eqnarray}}
\newcommand{\bars}{\begin{eqnarray*}}
\newcommand{\ears}{\end{eqnarray*}}
\newcommand{\nnu}{\nonumber \\}
\newtheorem{subn}{\name}
\renewcommand{\thesubn}{}
\newcommand{\bsn}[1]{\def\name{#1}\begin{subn}}
\newcommand{\esn}{\end{subn}}
\newtheorem{sub}{\name}[section]
\newcommand{\dn}[1]{\def\name{#1}}   
\newcommand{\bs}{\begin{sub}}
\newcommand{\es}{\end{sub}}
\newcommand{\bsl}[1]{\begin{sub}\label{#1}}
\newcommand{\bth}[1]{\def\name{Theorem}
\begin{sub}\label{t:#1}}
\newcommand{\blemma}[1]{\def\name{Lemma}
\begin{sub}\label{l:#1}}
\newcommand{\bcor}[1]{\def\name{Corollary}
\begin{sub}\label{c:#1}}
\newcommand{\bdef}[1]{\def\name{Definition}
\begin{sub}\label{d:#1}}
\newcommand{\bprop}[1]{\def\name{Proposition}
\begin{sub}\label{p:#1}}

\newcommand{\aand}{\quad\mbox{and}\quad}
\newcommand{\M}{{\cal M}}
\newcommand{\A}{{\cal A}}
\newcommand{\B}{{\cal B}}
\newcommand{\I}{{\cal I}}
\newcommand{\J}{{\cal J}}
\newcommand{\D}{\displaystyle}
\newcommand{\RR}{ I\!\!R}
\newcommand{\C}{\mathbb{C}}
\newcommand{\R}{\mathbb{R}}
\newcommand{\Z}{\mathbb{Z}}
\newcommand{\N}{\mathbb{N}}
\newcommand{\T}{{\rm T}^n}
\newcommand{\cuad}{{\sqcap\kern-.68em\sqcup}}
\newcommand{\abs}[1]{\mid #1 \mid}
\newcommand{\norm}[1]{\|#1\|}
\newcommand{\equ}[1]{(\ref{#1})}
\newcommand\rn{\mathbb{R}^N}
\renewcommand{\theequation}{\thesection.\arabic{equation}}
\newtheorem{definition}{Definition}[section]
\newtheorem{theorem}{Theorem}[section]
\newtheorem{proposition}{Proposition}[section]
\newtheorem{example}{Example}[section]
\newtheorem{proof}{proof}[section]
\newtheorem{lemma}{Lemma}[section]
\newtheorem{corollary}{Corollary}[section]
\newtheorem{remark}{Remark}[section]
\newcommand{\bremark}{\begin{remark} \em}
\newcommand{\eremark}{\end{remark} }
\newtheorem{claim}{Claim}


\newcommand{\rth}[1]{Theorem~\ref{t:#1}}
\newcommand{\rlemma}[1]{Lemma~\ref{l:#1}}
\newcommand{\rcor}[1]{Corollary~\ref{c:#1}}
\newcommand{\rdef}[1]{Definition~\ref{d:#1}}
\newcommand{\rprop}[1]{Proposition~\ref{p:#1}}
\newcommand{\BA}{\begin{array}}
\newcommand{\EA}{\end{array}}
\newcommand{\BAN}{\renewcommand{\arraystretch}{1.2}
\setlength{\arraycolsep}{2pt}\begin{array}}
\newcommand{\BAV}[2]{\renewcommand{\arraystretch}{#1}
\setlength{\arraycolsep}{#2}\begin{array}}
\newcommand{\BSA}{\begin{subarray}}
\newcommand{\ESA}{\end{subarray}}
\newcommand{\BAL}{\begin{aligned}}
\newcommand{\EAL}{\end{aligned}}
\newcommand{\BALG}{\begin{alignat}}
\newcommand{\EALG}{\end{alignat}}
\newcommand{\BALGN}{\begin{alignat*}}
\newcommand{\EALGN}{\end{alignat*}}
\newcommand{\note}[1]{\textit{#1.}\hspace{2mm}}
\newcommand{\Proof}{\note{Proof}}
\newcommand{\qeda}{\hspace{10mm}\hfill $\square$}
\newcommand{\qed}{\\
${}$ \hfill $\square$}
\newcommand{\Remark}{\note{Remark}}
\newcommand{\modin}{$\,$\\[-4mm] \indent}
\newcommand{\forevery}{\quad \forall}
\newcommand{\set}[1]{\{#1\}}
\newcommand{\setdef}[2]{\{\,#1:\,#2\,\}}
\newcommand{\setm}[2]{\{\,#1\mid #2\,\}}
\newcommand{\mt}{\mapsto}
\newcommand{\lra}{\longrightarrow}
\newcommand{\lla}{\longleftarrow}
\newcommand{\llra}{\longleftrightarrow}
\newcommand{\Lra}{\Longrightarrow}
\newcommand{\Lla}{\Longleftarrow}
\newcommand{\Llra}{\Longleftrightarrow}
\newcommand{\warrow}{\rightharpoonup}
\newcommand{
\paran}[1]{\left (#1 \right )}
\newcommand{\sqbr}[1]{\left [#1 \right ]}
\newcommand{\curlybr}[1]{\left \{#1 \right \}}
\newcommand{
\paranb}[1]{\big (#1 \big )}
\newcommand{\lsqbrb}[1]{\big [#1 \big ]}
\newcommand{\lcurlybrb}[1]{\big \{#1 \big \}}
\newcommand{\absb}[1]{\big |#1\big |}
\newcommand{\normb}[1]{\big \|#1\big \|}
\newcommand{
\paranB}[1]{\Big (#1 \Big )}
\newcommand{\absB}[1]{\Big |#1\Big |}
\newcommand{\normB}[1]{\Big \|#1\Big \|}
\newcommand{\produal}[1]{\langle #1 \rangle}

\newcommand{\thkl}{\rule[-.5mm]{.3mm}{3mm}}
\newcommand{\thknorm}[1]{\thkl #1 \thkl\,}
\newcommand{\trinorm}[1]{|\!|\!| #1 |\!|\!|\,}
\newcommand{\bang}[1]{\langle #1 \rangle}
\def\angb<#1>{\langle #1 \rangle}
\newcommand{\vstrut}[1]{\rule{0mm}{#1}}
\newcommand{\rec}[1]{\frac{1}{#1}}
\newcommand{\opname}[1]{\mbox{\rm #1}\,}
\newcommand{\supp}{\opname{supp}}
\newcommand{\dist}{\opname{dist}}
\newcommand{\myfrac}[2]{{\displaystyle \frac{#1}{#2} }}
\newcommand{\myint}[2]{{\displaystyle \int_{#1}^{#2}}}
\newcommand{\mysum}[2]{{\displaystyle \sum_{#1}^{#2}}}
\newcommand {\dint}{{\displaystyle \myint\!\!\myint}}
\newcommand{\q}{\quad}
\newcommand{\qq}{\qquad}
\newcommand{\hsp}[1]{\hspace{#1mm}}
\newcommand{\vsp}[1]{\vspace{#1mm}}
\newcommand{\ity}{\infty}
\newcommand{\prt}{\partial}
\newcommand{\sms}{\setminus}
\newcommand{\ems}{\emptyset}
\newcommand{\ti}{\times}
\newcommand{\pr}{^\prime}
\newcommand{\ppr}{^{\prime\prime}}
\newcommand{\tl}{\tilde}
\newcommand{\sbs}{\subset}
\newcommand{\sbeq}{\subseteq}
\newcommand{\nind}{\noindent}
\newcommand{\ind}{\indent}
\newcommand{\ovl}{\overline}
\newcommand{\unl}{\underline}
\newcommand{\nin}{\not\in}
\newcommand{\pfrac}[2]{\genfrac{(}{)}{}{}{#1}{#2}}

\def\ga{\alpha}     \def\gb{\beta}       \def\gg{\gamma}
\def\gc{\chi}       \def\gd{\delta}      \def\ge{\epsilon}
\def\gth{\theta}                         \def\vge{\varepsilon}
\def\gf{\phi}       \def\vgf{\varphi}    \def\gh{\eta}
\def\gi{\iota}      \def\gk{\kappa}      \def\gl{\lambda}
\def\gm{\mu}        \def\gn{\nu}         \def\gp{\pi}
\def\vgp{\varpi}    \def\gr{\rho}        \def\vgr{\varrho}
\def\gs{\sigma}     \def\vgs{\varsigma}  \def\gt{\tau}
\def\gu{\upsilon}   \def\gv{\vartheta}   \def\gw{\omega}
\def\gx{\xi}        \def\gy{\psi}        \def\gz{\zeta}
\def\Gg{\Gamma}     \def\Gd{\Delta}      \def\Gf{\Phi}
\def\Gth{\Theta}
\def\Gl{\Lambda}    \def\Gs{\Sigma}      \def\Gp{\Pi}
\def\Gw{\Omega}     \def\Gx{\Xi}         \def\Gy{\Psi}

\def\CS{{\mathcal S}}   \def\CM{{\mathcal M}}   \def\CN{{\mathcal N}}
\def\CR{{\mathcal R}}   \def\CO{{\mathcal O}}   \def\CP{{\mathcal P}}
\def\CA{{\mathcal A}}   \def\CB{{\mathcal B}}   \def\CC{{\mathcal C}}
\def\CD{{\mathcal D}}   \def\CE{{\mathcal E}}   \def\CF{{\mathcal F}}
\def\CG{{\mathcal G}}   \def\CH{{\mathcal H}}   \def\CI{{\mathcal I}}
\def\CJ{{\mathcal J}}   \def\CK{{\mathcal K}}   \def\CL{{\mathcal L}}
\def\CT{{\mathcal T}}   \def\CU{{\mathcal U}}   \def\CV{{\mathcal V}}
\def\CZ{{\mathcal Z}}   \def\CX{{\mathcal X}}   \def\CY{{\mathcal Y}}
\def\CW{{\mathcal W}} \def\CQ{{\mathcal Q}}
\def\BBA {\mathbb A}   \def\BBb {\mathbb B}    \def\BBC {\mathbb C}
\def\BBD {\mathbb D}   \def\BBE {\mathbb E}    \def\BBF {\mathbb F}
\def\BBG {\mathbb G}   \def\BBH {\mathbb H}    \def\BBI {\mathbb I}
\def\BBJ {\mathbb J}   \def\BBK {\mathbb K}    \def\BBL {\mathbb L}
\def\BBM {\mathbb M}   \def\BBN {\mathbb N}    \def\BBO {\mathbb O}
\def\BBP {\mathbb P}   \def\BBR {\mathbb R}    \def\BBS {\mathbb S}
\def\BBT {\mathbb T}   \def\BBU {\mathbb U}    \def\BBV {\mathbb V}
\def\BBW {\mathbb W}   \def\BBX {\mathbb X}    \def\BBY {\mathbb Y}
\def\BBZ {\mathbb Z}

\def\GTA {\mathfrak A}   \def\GTB {\mathfrak B}    \def\GTC {\mathfrak C}
\def\GTD {\mathfrak D}   \def\GTE {\mathfrak E}    \def\GTF {\mathfrak F}
\def\GTG {\mathfrak G}   \def\GTH {\mathfrak H}    \def\GTI {\mathfrak I}
\def\GTJ {\mathfrak J}   \def\GTK {\mathfrak K}    \def\GTL {\mathfrak L}
\def\GTM {\mathfrak M}   \def\GTN {\mathfrak N}    \def\GTO {\mathfrak O}
\def\GTP {\mathfrak P}   \def\GTR {\mathfrak R}    \def\GTS {\mathfrak S}
\def\GTT {\mathfrak T}   \def\GTU {\mathfrak U}    \def\GTV {\mathfrak V}
\def\GTW {\mathfrak W}   \def\GTX {\mathfrak X}    \def\GTY {\mathfrak Y}
\def\GTZ {\mathfrak Z}   \def\GTQ {\mathfrak Q}

\font\Sym= msam10 
\newcommand{\bdw}{\prt\Gw\xspace}

 \centerline{\large \bf Liouville theorems for elliptic equations involving regional}
 
 \smallskip
 
 \centerline{\large \bf fractional Laplacian  with order in $(0,\,1/2]$}
 
 \bigskip\medskip

\centerline{   Huyuan Chen \quad {\rm and}\quad Yuanhong Wei }
\medskip
{\footnotesize
} 

\begin{abstract}
In this paper, some elliptic equation in a bounded open domain in $\mathbb{R}^N$ ($N\geq 2$) with $C^2$ boundary $\partial\Omega$ is considered. The problem is driven by the regional fractional Laplacian, the  infinitesimal generator  of the censored symmetric $2\alpha$-stable process in $\Omega$.  Probability theory asserts that the censored $2\alpha$-stable process can not approach the boundary when $\alpha\in(0,\frac12]$.
 For $\alpha\in (0,\frac12]$,   our purpose in this article is to show that  non-existence of  solutions bounded from above or bounded from below for the particular  Poisson problem
$$
   (-\Delta)^\alpha_\Omega   u= 1 \quad   {\rm in}\ \, \,  \Omega 
 $$
and
non-existence of nonnegative nontrivial solutions  of the Lane-Emden equation
$$
\arraycolsep=1pt\left\{
\begin{array}{lll}
 \displaystyle  (-\Delta)^\alpha_\Omega   u=u^p\quad & {\rm in}\ \,  \Omega,\\[1.5mm]
\phantom{ (-\Delta)^\alpha  }
 \displaystyle   u=0\quad & {\rm on}\   \partial\Omega.
 \end{array}\right.
 $$

\end{abstract}

\vspace{1mm}

\vspace{1mm}
\hspace{.05in}
\medskip

\setcounter{equation}{0}
\section{Introduction}

In the last few years,  there has been an increasing interest in the study of nonlocal problems which are now widely used in physics, operation research, queuing theory, mathematical finance, risk estimation, and others (see \cite{B,CSS,GM,S}).  Given a  bounded open Lipchitz set $\Omega\subset\mathbb{R}^N$ with $N\geq 2$,   a symmetric $2\alpha$-stable process $X=\{X_t\}$ in $\mathbb{R}^N$  killed upon leaving the set $\Omega$ is called a symmetric $2\alpha$-stable process  in $\Omega$, denote it by $X^\Omega$,  see \cite{CKS}.  The Dirichlet form of $X^\Omega$ on $L^2(\Omega,dx)$ is $(\CC,\CF^\Omega)$, where
$$\CF^\Omega=\Big\{u\in L^2(\R^N):\, u=0\ \ {\rm q.e.\ on}\ \mathbb{R}^N\setminus\Omega,\ \, \int_{\mathbb{R}^N}\int_{\mathbb{R}^N}\frac{(u(x)-u(y))^2}{|x-y|^{N+2\alpha}}dxdy<+\infty\Big\}$$
and
$$\CC(u,v)=\frac12 c_{N,\alpha} \int_{\Omega}\int_{\Omega}\frac{(u(x)-u(y))(v(x)-v(y)) }{|x-y|^{N+2\alpha}}dxdy+\int_\Omega \kappa_\alpha(x) u(x)v(x)dx,\ \forall\, u,v\in  \CF^\Omega.$$
Here q.e. is the abbreviation for quasi-everywhere (cf. \cite{BBC}, \cite{FO}), and  the density of the killing measure of $X^\Omega$ $$\displaystyle \kappa_\alpha(x)=c_{N,\alpha} \int_{\mathbb{R}^N\setminus \Omega} |x-y|^{-N-2\alpha}dy.$$ 
Some properties of killing measure $\kappa_\alpha$ are given in Appendix.   

Elliptic differential equations and diffusion processes play significant roles in the theory of partial differential equations and in probability theory, respectively. There are close relationships between these two subjects. Sometimes for an elliptic operator, there is a diffusion process associated with it, so that the elliptic operator is the infinitesimal generator of the diffusion process.

The inhomogenous elliptic problem associated to this kind of process is 
 \begin{equation}\label{eq 1.1-0}
\arraycolsep=1pt\left\{
\begin{array}{lll}
 \displaystyle  (-\Delta)^\alpha    u=f\quad & {\rm in}\ \,  \Omega,\\[1mm]
\phantom{ (-\Delta)^\alpha  }
 \displaystyle   u=0\quad & {\rm in}\   \mathbb{R}^N\setminus \Omega,
 \end{array}\right.
\end{equation}
where
$$ (-\Delta)^\alpha   u(x)=c_{N,\alpha}{\rm p.v.}\int_{\mathbb{R}^N}\frac{ u(x)-
u(z)}{|z-x|^{N+2\alpha}}\,dz$$
is the the fractional Laplacian operator and $c_{N,\alpha}>0$ is the normalized constant, see  \cite{EGE}.

Define $\tau_\Omega=\inf\{t>0: X_t\not\in\Omega\}$ and note that $\lim_{t\to\tau_\Omega}X_t$ exists and typically belongs to $\Omega$. Bogdan, Burdzy and Chen \cite{BBC} (also see Guan and Ma \cite{GM}) extended $X^\Omega$ beyond its lifetime $\tau_\Omega$
and obtained a version of the strong Markov process, named censored symmetric stable process. A censored stable process in an open set $\Omega\subset\mathbb{R}^N$ is obtained from the symmetric stable process by suppressing its jumps from $\Omega$ to $\mathbb{R}^N\backslash\Omega$. In fact, the censored stable process is a stable process forced to stay inside $\Omega$ and the Dirichlet form has no killing term.  
It is known that the censored stable process has the generator
the regional fractional Laplacian defined in $\Omega$ defined as
$$ (-\Delta)^\alpha_\Omega  u(x)=\lim_{\varepsilon\to0^+} (-\Delta)_{\Omega,\varepsilon}^\alpha u(x)$$
with
$$
(-\Delta)_{\Omega,\varepsilon}^\alpha  u(x)=c_{N,\alpha}\int_{\Omega\setminus B_\varepsilon(x)}\frac{ u(x)-
u(z)}{|z-x|^{N+2\alpha}}\,dz.
$$

It is believed that censored stable processes deserve to be studied because the classical counterpart, killed Brownian motion, is an important model in both pure mathematics and in applied probability (see \cite{BBC}). When $\alpha\in(0,\, 1)$, Kulczycki \cite{K} establishes the upper and the lower bound estimates of Green function for a bounded open domain. When $\alpha\in(\frac12,\, 1)$, \cite{CK,CKS} give estimates on the heat kernel and Green kernel related to
the   regional fractional Laplacian. \cite{G} provides a version of integration by formula for regional fractional  Laplacian and \cite{C} extends this formula to solve regional fractional  problem with kinds of inhomogenous terms. More study of regional fractional laplacian with $\alpha\in (\frac12,1)$ could see  \cite{CH} boundary blowing-up solutions, \cite{F} boundary regularity, \cite{T} weak solution for Poisson problem, \cite{D} involving Hardy potential and more related properties \cite{DFW,BGU, AT,TC}.

It should be pointed out that the censored stable process
has some interesting properties, which imply some differences between  $(\frac12,\, 1)$ and $(0, \, \frac12]$  (see \cite[Theorem 1]{BBC}):
\begin{itemize}
\item[(i)\ ]  if $\alpha\in(\frac12,\, 1)$, the censored symmetric $2\alpha$-stable process in $\Omega$ has a finite lifetime and will approach $\partial \Omega$;
\item[(ii)]  if $\alpha\in(0, \, \frac12]$, the censored symmetric $2\alpha$-stable process in $\Omega$   is conservative and will never approach  $\partial \Omega$.
\end{itemize}\smallskip
This probability phenomena indicates that the Poisson problem 
\begin{equation}\label{eq 1.1-y}
\arraycolsep=1pt\left\{
\begin{array}{lll}
 \displaystyle  (-\Delta)^\alpha_\Omega   u=f\quad & {\rm in}\ \,  \Omega,\\[1.5mm]
\phantom{ (-\Delta)^\alpha  }
 \displaystyle   u=0\quad & {\rm on}\   \partial\Omega
 \end{array}\right.
\end{equation}
would have very different solution's structures between $\alpha\in(\frac12,1)$ and $\alpha\leq (0,\frac12]$. Moreover,  \cite{M} shows that the space $H^\alpha_0(\Omega)$ has zero trace when $\alpha\in(\frac12,\, 1)$ and  has no zero trace when $\alpha\in (0,\, \frac12]$; \cite{BD} shows that the optimal constant $\CA_{N,\alpha}$ in
the fractional Hardy inequality
$$\int_{\BBR^N_+}\int_{\BBR^N_+}\frac{(u(x)-u(y))^2 }{|x-y|^{N+2\alpha}}dxdy\geq \CA_{N,\alpha} \int_{\BBR^N_+} x_N^{-2\alpha} u^2(x) dx\quad{\rm for\ any} \ u\in C_c(\BBR^N_+),$$
is positive if $\alpha\in(\frac12, \,1)$ and zero if $\alpha\in (0,\, \frac12]$,
where $\mathbb{R}^N_+=\{x=(x',x_N)\in \mathbb{R}^{N-1}\times \mathbb{R}: x_N>0\}$. 
This property produces  a big challenge to consider the solution of the Poisson problem (\ref{eq 1.1-y}) for $\alpha\in(0,\, \frac12]$.

Via building the formula of integral by part  and related  embedding results,   \cite{C} obtains the existence of  solutions to (\ref{eq 1.1-y}) for $\alpha\in(\frac12,1)$ when $\Omega$ is a bounded regular domain. Precisely, we have that 

\begin{proposition}\label{pr 1}\cite{C}
 Assume that $\alpha\in(\frac12,1)$,  $\Omega$ is a bounded $C^2$ domain in $\R^N$ with $N\geq 2$ and $h=0$ on $\partial\Omega$.

 $(i)$ If $f$ is H\"older continuous, then  (\ref{eq 1.1-y}) has a classical solution $u_f\in C_0(\Omega)$;

 $(ii)$ If $f\in L^2(\Omega)$, then (\ref{eq 1.1-y}) has a weak solution $u_f\in H^\alpha_0(\Omega)$,
 where $H^\alpha_0(\Omega)$ is the Hilbert space, which is the closure of $C_c^\infty(\Omega)$ under the norm
 $$\norm{u}_{\alpha}:= \sqrt{\int_{\Omega}\int_{\Omega}\frac{(u(x)-u(y))^2 }{|x-y|^{N+2\alpha}}dxdy}.$$
\end{proposition}
In spite of some known results, the theory of elliptic problems driven by regional fractional Laplacian is far from mature and adequate. 

 For $\alpha\in(0,\frac12]$, inversely, our interest in this article is to show the nonexistence of solution to particular problem 
\begin{equation}\label{eq 1.1}
  (-\Delta)^\alpha_\Omega   u=1\quad  {\rm in}\ \,  \Omega.
\end{equation}
Here we don't restrict our problem with a Dirichlet boundary condition, so we may consider solutions of (\ref{eq 1.1}) with a wild behavior at the boundary.  Here a viscosity solution $u\in C(\Omega)$ of (\ref{eq 1.1})
is bounded from above (\,bounded from below resp.\,)  if
$$\sup_{x\in\Omega} u(x)<+\infty\quad (\, \inf_{x\in\Omega} u(x)>-\infty\, ).$$

Our main results on Poisson problem (\ref{eq 1.1})
with $\alpha\in(0,\, \frac12]$ states as follows.

\begin{theorem}\label{teo 1}
 Assume that the integer $N\geq1$, $\alpha\in(0, \frac12]$, $\Omega$ is a bounded   domain in $\BBR^N$, which is $C^2$ if $N\geq 2$ and an open interval if $N=1$.

Then   problem (\ref{eq 1.1}) has no  solution bounded from above or bounded from below.

\end{theorem}

 Our second purpose is to show the Liouville type theorem for 
  the Lane-Emden equation with the regional fractional Laplacian
 \begin{equation}\label{eq 1.2}
\arraycolsep=1pt\left\{
\begin{array}{lll}
 \displaystyle  (-\Delta)^\alpha_\Omega   u=u^p\quad & {\rm in}\ \,  \Omega,\\[1.5mm]
\phantom{ (-\Delta)^\alpha  }
 \displaystyle   u=0\quad & {\rm on}\   \partial\Omega,
 \end{array}\right.
\end{equation}
where $p>0$. 

\begin{theorem}\label{teo 2}
Assume that the integer $N\geq1$, $\alpha\in(0, \frac12]$, $\Omega$ is a bounded   domain in $\BBR^N$, which is $C^2$ if $N\geq 2$ and an open interval if $N=1$.
Then problem (\ref{eq 1.2})  has no nonnegative  nontrivial  solution for any  $p>0$.
\end{theorem}

It is worth noting that equation (\ref{eq 1.1}) replacing the regional fractional Laplacian by Laplacian or fractional Laplacian has a unique solution under the Dirichlet boundary condition; 
 equation (\ref{eq 1.2}) replacing the regional fractional Laplacian by Laplacian or fractional Laplacian could have a classical solution if $p\in(1,\frac{N+2\alpha}{N-2\alpha})$ Sobolev subcritical with $s\in(0,1]$. 

Our results show a very different phenomena  of nonexistence,   is to investigate non-existence of viscosity super-solutions, which, to some extent, can be seen as an interpretation of the result in \cite{BBC}  and an counterpart of \cite{C} for $\alpha\in(0,\, \frac12]$. 

Our idea to prove the Liouville theorem is based on the non-existence boundary blowing up solution of 
 \begin{equation}\label{eq 1.3}
\arraycolsep=1pt\left\{
\begin{array}{lll}
 \displaystyle  (-\Delta)^\alpha_\Omega   u=0\quad & {\rm in}\ \,  \Omega,\\[1.5mm]
\phantom{ (-\Delta)^\alpha  }
 \displaystyle   u=+\infty\quad & {\rm on}\   \partial\Omega.
 \end{array}\right.
\end{equation}
It is worth noting that 
\begin{itemize}
\item[(i)\ ]  when $s=1$,  there is no  the classical harmonic functions with boundary blowing up;
\item[(ii)]  For $\alpha\in (0,1)$, there is harmonic functions with boundary blowing-up 
\begin{equation}\label{eq 1.3}
\arraycolsep=1pt\left\{
\begin{array}{lll}
 \displaystyle \quad\  (-\Delta)^\alpha    u=0\quad & {\rm in}\ \,  \Omega,\\[1.5mm]
\phantom{   }
 \displaystyle  \lim_{x\in\Omega, x\to\partial\Omega} u=+\infty\quad & {\rm on}\   \partial\Omega,\\[1.5mm]
\phantom{ (-\Delta)^\alpha -\ \  }
 \displaystyle
 u=0\quad &{\rm in}\ \R^N\setminus \Omega.
 \end{array}\right.
\end{equation}
which could see \cite{CFQ,A}, the solution has the behavior $\rho^{s-1}(x)$, where 
$\rho(x)={\rm dist}(x,\partial\Omega)$.
\end{itemize}

The key point is to obtain  estimates  of $(-\Delta)^\alpha_\Omega V_\tau$ for $\alpha\in(0,\frac12]$ near the boundary,
 where
\begin{equation}\label{1.1}
V_\tau(x)=\left\{ \arraycolsep=1pt
\begin{array}{lll}
 \rho(x)^{\tau}\ \ \ &\text{ for } \, x\in A_\delta, \\[2mm]
 l(x)\ \ \  &\text{ for } \, x\in \Omega\setminus A_\delta
 \end{array}
\right.
\end{equation}
and
\begin{equation}\label{1.1-*}
V_0^*(x)=\left\{ \arraycolsep=1pt
\begin{array}{lll}
 -\ln \rho(x)\ \ \ &\text{ for } \, x\in A_\delta, \\[2mm]
 l(x)\ \ \  &\text{ for } \, x\in \Omega\setminus A_\delta
 \end{array}
\right.
\end{equation}
with $A_\delta=\{x\in \Omega: \rho(x)<\delta\}$,  $\delta\leq \frac12$,    $\tau\in (-1,2\alpha)$ and the function $l$  being  positive  such that
$V_\tau$ is $C^2$ in $\Omega$. It is worth noting that when $\tau\in (2\alpha-1, 0)$ with $\alpha\in(0,\frac12)$, $(-\Delta)^\alpha_\Omega V_\tau$
blows up near the boundary positively, which leads to the non-existence of Poisson problem (\ref{1.1}). \smallskip

The rest of this paper is organized as follows. In Section 2, we recall the definition of viscosity solution and do the estimates of
$(-\Delta)^\alpha_\Omega V_\tau$. In Section 3,  we prove the non-existence of Poisson problem (\ref{eq 1.1}) by contradiction
that a boundary blow-up super solution could be constructed if   (\ref{eq 1.1}) has a super solution.
  Finally, we annex properties of killing measure $\kappa_\alpha$.




\mysection{Preliminary}

\subsection{ Viscosity solution }

We start with the  definition  of viscosity solution, inspired by the definition of
viscosity sense for nonlocal problem in \cite{CS}.

\begin{definition}\label{de 2.2}
We say that a  function $u\in C( \Omega)$ is a viscosity super solution (sub-solution)
 of
 \begin{equation}\label{eq 2.1}
 \left\{ \arraycolsep=1pt
\begin{array}{lll}
 \displaystyle  (-\Delta)^{\alpha}_\Omega u=f\quad & {\rm in}\ \,   \Omega,\\[1.5mm]
\phantom{ (-\Delta)^\alpha  }
 \displaystyle   u=h\quad & {\rm on}\    \partial  \Omega,
\end{array}\right.
\end{equation}
if $u\geq h$ (resp. $u\leq h$) on $\partial  \Omega$ and
for every point $x_0\in\Omega$ and some  neighborhood $V$ of
$x_0$ with $\bar V\subset \Omega $ and for any $\varphi \in
C^2(\bar V)$ such that $u(x_0)=\varphi(x_0)$ and $x_0$ is the minimum (resp. maximum) point of $u-\varphi$ in $V$,  let
\begin{eqnarray*}
\tilde u  =\left\{ \arraycolsep=1pt
\begin{array}{lll}
\varphi\ \ \ & \rm{in}\ \, &V,\\[1mm]
u \ \ & \rm{in}\ \, &\Omega\setminus V,
\end{array}
\right.
\end{eqnarray*}
we have that
$$(-\Delta)^{\alpha }_\Omega\tilde u(x_0) \geq f(x_0)\quad (resp.  \ (-\Delta)^{\alpha }_\Omega\tilde u(x_0) \le f(x_0)).$$

We say that $u$ is a
viscosity solution of (\ref{eq 2.1}) if it is  a viscosity super-solution and also a viscosity sub-solution of (\ref{eq 2.1}).
\end{definition}

\begin{theorem}\label{comparison}
Assume that the functions  $f:\Omega\to\BBR$, $h:\partial\Omega\to\BBR$ are
continuous. Let $u$ and $v$ be a viscosity super-solution and
sub-solution  of (\ref{eq 2.1}) respectively. Then
\begin{equation}\label{2.1}
 v  \le u \quad{\rm in}\quad  \Omega.
\end{equation}

\end{theorem}
\noindent{\bf Proof.}  Let us define $w=u-v$, then
\begin{equation}\label{eq 2.2}
 \left\{ \arraycolsep=1pt
\begin{array}{lll}
 \displaystyle  (-\Delta)^{\alpha}_\Omega w \geq 0\quad\, & {\rm in}\ \,  \Omega,\\[1.5mm]
\phantom{ (-\Delta)^\alpha   }
 \displaystyle  w\geq 0\quad\, & {\rm on}\   \partial  \Omega.
\end{array}\right.
\end{equation}
If (\ref{2.1}) fails, then there exists $x_0\in\Omega$ such that
$$w(x_0)=u(x_0)-v(x_0)=\min_{x\in\Omega}w(x)<0,$$
  then in the viscosity sense,
\begin{equation}\label{y 2.1}
 (-\Delta)^{\alpha}_\Omega w(x_0)\geq 0.
\end{equation}
 Since $w$ is a viscosity super solution  $x_0$ is the minimum point in $\Omega$ and $w\ge 0$ on $\partial\Omega$, then we can
take a small neighborhood $V_0$ of $x_0$ such that $\tilde w=w(x_0)$ in $V_0$,
From (\ref{y 2.1}), we have that
$$ (-\Delta)^{\alpha}_\Omega \tilde  w(x_0)\geq 0.$$
But the definition of regional fractional Laplacian implies that
$$(-\Delta)^{\alpha}_\Omega \tilde w(x_0)=   \int_{\Omega\setminus V_0} \frac{w(x_0)-w(y)}{|x_0-y|^{N+2\alpha}}dy<0,$$
which is impossible.  \qquad $\Box$

\begin{remark}\label{rm 2.1-1}
Let $u$ be a continuous function in $\Omega$ and  $x_0$ is a minimum point of $u$, then
$(-\Delta)^\alpha_\Omega u(x_0)>0$ in the viscosity sense.

\end{remark}

\subsection{ Estimates for boundary blowing up functions }  It is known the derivatives on  the boundary blowing-up functions is related the one-dimensional problem and we start this subsection by the analysis of one-dimensional regional fractional problem on the half-line. Let
\begin{equation}\label{3.1}
  w_\tau(t)=  t^{\tau} \quad {\rm and}\quad w_0^*(t)=-\ln  t\quad{\rm for }\ \, t>0,
\end{equation}
where $\tau>-1$.
\begin{lemma}\label{lm 3.1}
Assume that $\alpha\in(0,\, 1)$,  $\tau\in(-1,\, 2\alpha)$ and $w_\tau, w_0^*$ are defined in (\ref{3.1}). Let
$$I_\alpha=\left\{ \arraycolsep=1pt
\begin{array}{lll}
  (0,\, 2\alpha-1)\quad & {\rm if}\ \alpha\in(\frac12,\, 1),\\[1.5mm]
\phantom{     }
 \displaystyle  (2\alpha-1,\, 0)\quad &   {\rm if}\ \alpha\in(0,\, \frac12).
\end{array}\right.
$$
Then
\begin{equation}\label{2.11}
(-\Delta)^{\frac12}_{\R_+} w_0^*=0\quad \text{ for }\, t>0
\end{equation}
and
$$(-\Delta)^{\alpha}_{\R_+} w_\tau=c_\alpha(\tau) t^{\tau-2\alpha}\quad \text{ for }\, t\in\R_+,$$
where
\begin{equation}\label{3.2}
c_\alpha(\tau)
\left\{ \arraycolsep=1pt
\begin{array}{lll}
 >0 \qquad{\rm for}\ \ \tau\in I_\alpha,\\[1.5mm]
=0   \qquad{\rm for}\ \ \tau=2\alpha-1,\\[1.5mm]
<0 \qquad{\rm for}\ \ \tau\in (-1,2\alpha)\setminus \bar I_\alpha.
\end{array}\right.
\end{equation}
\end{lemma}
\noindent {\bf Proof.}  {\it Estimates of $c_\alpha(\tau)$:} For $\epsilon\in(0,\frac18)$, by change variable we have that
\begin{eqnarray*}
(-\Delta)_{\R_+}^\alpha  w_\tau(t) &=&c_{1,\alpha}t^{\tau-2\alpha}\lim_{\epsilon\to0^+}\int_{(0,\infty)\setminus (1-\epsilon,1+\epsilon)}
\frac{ 1-s^\tau}{|1-s|^{1+2\alpha}}\,ds\\
&=&c_{1,\alpha}t^{\tau-2\alpha}\lim_{\epsilon\to0^+}\Big(\int_0^{1-\epsilon}\frac{1-s^{\tau}}{(1-s)^{1+2\alpha}} ds-\int_{ 1+\epsilon}^{\infty}\frac{1-s^{-\tau}}{(1-s^{-1})^{1+2\alpha}}s^{-2\alpha-1+\tau} ds \Big) \\
&=& c_{1,\alpha}t^{\tau-2\alpha}\lim_{\epsilon\to0^+}\Big(\int_0^{1-\epsilon}\frac{(1-s^{\tau})(1-s^{2\alpha-\tau-1})}{(1-s)^{1+2\alpha}} ds+\int_{1-\epsilon}^{\frac1{1+\epsilon}}\frac{(1-s^{\tau})s^{2\alpha-1-\tau}}{(1-s)^{1+2\alpha}} ds \Big),
\end{eqnarray*}
where $1-\epsilon<\frac1{1+\epsilon}<1-\epsilon+\epsilon^2$ and
\begin{eqnarray*}
\Big| \int_{1-\epsilon}^{\frac1{1+\epsilon}}\frac{(1-s^{\tau})s^{2\alpha-1-\tau}}{(1-s)^{1+2\alpha}} ds  \Big|
 &\leq&  c_1|\tau| \int_{1-\epsilon}^{\frac1{1+\epsilon}} (1-s)^{-2\alpha} ds
 \\&\leq& c_2 \epsilon^{2-2\alpha}\to0\quad{\rm as}\ \, \epsilon\to0^+.
 \end{eqnarray*}

As a consequence,
$$(-\Delta)_{\R_+}^\alpha  w_\tau(t)=c_{1,\alpha}\gamma(\alpha,\tau) t^{\tau-2\alpha},  $$
where
$$\gamma(\alpha,\tau)=\int_0^1\frac{(1-s^\tau)(1-s^{2\alpha-1-\tau})}{(1-s)^{1+2\alpha}}ds.$$
Note that
\begin{eqnarray*}
 &\gamma(\alpha,\tau)>0\ \ &\text{ if and only if }\ \ \tau(2\alpha-1-\tau)>0,\\
 &\gamma(\alpha,\tau)< 0\ \ &\text{ if and only if }\ \ \tau(2\alpha-1-\tau)< 0,\\
 &\gamma(\alpha,\tau)= 0\  \ &\text{ if   }\quad \ \tau=0\ \ {\rm or}\ \ \tau=2\alpha-1.\quad
 \end{eqnarray*}
Therefore, we conclude   (\ref{3.2}) with $c_\alpha(\tau)=c_{1,\alpha}\gamma(\alpha,\tau)$.\smallskip

{\it Proof of  $(-\Delta)_{\R_+}^{1/2} w_0^*=0$. }  By change variable we have that
\begin{eqnarray*}
(-\Delta)_{\R_+}^{1/2}  w_0^*(t) &=&c_{1,1/2}t^{-1}\lim_{\epsilon\to0^+}\int_{(0,\infty)\setminus (1-\epsilon,1+\epsilon)}
\frac{ \ln s }{|1-s|^2}\,ds\\
&=&c_{1,1/2}t^{-1}\lim_{\epsilon\to0^+}\Big(\int_0^{1-\epsilon}\frac{\ln s }{(1-s)^2} ds+\int_0^{\frac1{1+\epsilon}} \frac{-\ln s }{(1-s)^2}  ds \Big) \\
&=& c_{1,1/2}t^{-1} \lim_{\epsilon\to0^+} \int_{1-\epsilon}^{\frac1{1+\epsilon}}\frac{-\ln s }{(1-s)^2}  ds,
\end{eqnarray*}
where $1-\epsilon<\frac1{1+\epsilon}<1-\epsilon+\epsilon^2$ and for $c_3,c_4>0$,
\begin{eqnarray*}
\Big| \int_{1-\epsilon}^{\frac1{1+\epsilon}}\frac{-\ln s }{(1-s)^2}  ds  \Big|
 &\leq&  c_3  \int_{1-\epsilon}^{\frac1{1+\epsilon}} (1-s)^{-1} ds\\&\leq& c_4 \epsilon \to0\quad{\rm as}\ \, \epsilon\to0^+.
 \end{eqnarray*}
 Therefore,  $(-\Delta)^{\frac12}_{\R_+}w_0^*=0$ in $\R_+$. We complete the proof.
\hfill$\Box$

     \medskip

\begin{proposition}\label{pr 1.1}
 Assume that  $\alpha\in(0, \, 1)$, $I_\alpha=(-1,\, 2\alpha)$, $\Omega$ is a bounded $C^2$ domain in $\R^N$ with $N\geq 1$
 and $V_\tau$ is given in (\ref{1.1}) with $\tau\in I_\alpha$. Let
\begin{equation}\label{interval}
I_\alpha^*=\left\{ \arraycolsep=1pt
\begin{array}{lll}
(2\alpha-1,\, 0)  \quad&{\rm if}\ \ \alpha\in(0,\,\frac12),\\[1.5mm]
(0,\, 2\alpha-1)   \quad&{\rm if}\ \  \alpha\in (\frac12,\, 1),\\ [1.5mm]
\emptyset \quad&{\rm if}\ \  \alpha=\frac12
\end{array}\right.\quad{\rm and}\quad  I_\alpha^\#=\left\{ \arraycolsep=1pt
\begin{array}{lll}
\bar I_\alpha^* \quad&{\rm if}\ \ \alpha\not=\frac12,\\[1.5mm]
\{0\}   \quad&{\rm if}\ \  \alpha =\frac12.
\end{array}\right.
\end{equation}

 Then
 there exist $\delta_1\in (0,\delta)$ and  constant $c_5>1$ such that:
\\
 $(i)$  If $\tau\in I_\alpha^*$, then
$$
\frac1{c_5}\rho(x)^{\tau-2\alpha }\leq (-\Delta)^{\alpha}_\Omega V_\tau(x)\leq
c_5\rho(x)^{\tau-2\alpha }\ \ \mbox{for all}\,\, x\in A_{\delta_1}.$$
$(ii)$\ If
 $\tau\in I_\alpha\setminus  I_\alpha^\#$, then
$$
\frac1{c_5}\rho (x)^{\tau-2\alpha }\leq -(-\Delta)^{\alpha}_\Omega V_\tau(x)\leq
c_5\rho(x)^{\tau-2\alpha }\ \ \mbox{for all}\,\, x\in A_{\delta_1}.$$
 $(iii)$ \ If
 $\tau=2\alpha-1$ and $\alpha\not=\frac12$, then
$$|(-\Delta)^{\alpha}_\Omega V_\tau(x)| \leq
c_5 (\rho(x)^{\tau-2\alpha+1 }+ \rho(x)^\tau+1)\ \ \mbox{for all}\,\, x\in A_{\delta_1}.$$

\end{proposition}
 \noindent{\bf  Proof.}
 When $N=1$, we assume that $(0,\delta_0)\subset \Omega$,  for $x_1\in(0,\frac{\sigma_0}4)$, we have that
 \begin{eqnarray*}
(-\Delta)_\Omega^\alpha \nonumber V_\tau(x)&=&  c_{1,\alpha} {\rm p.v.} \int_0^{\delta_0}
 \frac{x^{\tau}-y^{\tau}}{|x-y|^{1+2\alpha}}  dy +\int_{\Omega \setminus (0,\delta_0)}
 \frac{x^{\tau}-y^{\tau}}{|x-y|^{1+2\alpha}}  dy
\\&=&   c_{1,\alpha} x^{\tau-2\alpha}\int^{\frac{\delta_0}x}_0
 \frac{1-t^{\tau}}{|1-t|^{1+2\alpha}}  dt+O(1)(x^\tau+1)
 \\&=&x^{\tau-2\alpha} \left(c_\alpha(\tau)-(\frac{\delta_0}{x})^{\tau-2\alpha}\right)+O(1)(x^\tau+1)
 \\&=& c_\alpha(\tau)x^{\tau-2\alpha}+O(1)(x^\tau+1).
\end{eqnarray*}

 Now we  concentrate on the case that $N\geq 2$.
  By compactness we prove that the corresponding inequality holds in a neighborhood of any point $\bar x\in\partial\Omega$ and  without loss of generality we may assume that $\bar
x=0$ and we further assume that  $-e_1$ is the outer normal
vector of $\Omega$ at $0$, where $e_1=(1,0,\cdots)\in\R^N$.

 For  a given $0<\eta\le \delta$, we define
$$
Q_\eta=\{y=(y_1,y')\in\R\times\R^{N-1},\, 0<y_1<\eta,\, |y'|<\eta\}$$
and
$$
\tilde Q_\eta=\{y=(y_1,y')\in\R\times\R^{N-1},\, |y_1|<\eta,\, |y'|<\eta\}.$$

Since $\partial \Omega$ is $C^2$, then there is a $C^2$ function $\varphi:\R^{N-1}\to\R$  such that
  $(z_1, z')\in \Omega\cap \tilde Q_{2\eta}$ if and only if $z_1=\varphi(z')$ for $|z'|<2\eta$ for $\eta>0$ small enough. In fact, we take $\eta>0$ small enough, then for any $z\in A_t:=\{x\in \Omega:\, \rho(x)=t\}$ with $t\in(0,\eta]$, there exists a unique point $z_0\in \partial \Omega$ such that $|z-z_0|=t$, $z-z_0$ is a normal vector at $z_0$
and  $A_t$ is $C^2$.

Let $\Phi:\R^N\to \R^N$ be a $C^2$ diffeomorphism such that
$$\Phi(z)=(z', z_1+\varphi(z'))\quad{\rm for}\ \, z=(z_1,\, z')\in   Q_\eta.$$
Note that
$$\Phi(te_1)=te_1\ \ {\rm for}\ \, t\in(0,\eta), \quad \quad \Phi(0,y')\in \partial\Omega\cap \tilde Q_\eta\quad{\rm for}\ \, |y'|<\eta$$
and
$$\tilde Q_{\frac{\eta}2}\cap \Omega\subset  \Phi(Q_\eta) \subset \tilde Q_{2\eta}\cap \Omega.$$
Moreover,  $\rho(\Phi(y))=y_1+O(|y'|^2)$ for $y\in Q_\eta$ and there exists $c_6>0$ such that
$$\Phi(0)=0,\ \ |D\Phi(0)-E_1|\leq c_6 |z'|\quad{\rm and}\quad |\Phi(y)-y|\leq c_6|y'|^2\quad{\rm for }\ y=(y_1,y')\in Q_\eta,$$
where $E_1=(d_{ij})_{N\times N}$ is the unit matrix, i.e. $d_{ii}=1$ and $d_{ij}=0$ if $i\not=j$.

Let $x_1\in(0,\eta/4)$,
$x=(x_1,0)$   be  a generic point in $A_{\eta/4}$. We observe that $|x-\bar
x|=\rho(x)=x_1$.
By definition we have
\begin{eqnarray}
(-\Delta)^\alpha V_\tau(x)=c_{N,\alpha}{\rm p.v.} \int_{\Phi(Q_{ \eta})}\frac{x_1^{\tau}-\rho^\tau(z) }{|x-z|^{N+2\alpha}}dz+c_{N,\alpha} \int_{\Omega\setminus
\Phi(Q_{ \eta})}\frac{V_{\tau}(x)-V_\tau(z) }{|x-z|^{N+2\alpha}}dz \label{cotadelta}
\end{eqnarray}
and we see that
\begin{equation}\label{E1(x1)+}
|\int_{\Omega\setminus\Phi(Q_{ \eta})}\frac{V_{\tau}(x)-V_\tau(z) }{|x-z|^{N+2\alpha}}dz|
\le c_7(x_1^\tau+1),
\end{equation}
where the constant $c_7$ is independent of $x$. Thus we only need to study the asymptotic behavior of the first integral, that from now on, we denote by
$$\BBE(x_1)={\rm p.v.} \int_{\Phi(Q_{ \eta})}\frac{x_1^{\tau}-\rho^\tau(z) }{|x-z|^{N+2\alpha}}dz.$$
Let $z=\Phi(y)$, then we have that
$$\BBE(x_1)={\rm p.v.} \int_{ Q_{ \eta}}\frac{x_1^{\tau}-y_1^{\tau} }{|x-\Phi(y)|^{N+2\alpha}}
|D\Phi(y)| dy+\int_{ Q_{ \eta}}\frac{O(|y'|^2) }{|x-\Phi(y)|^{N+2\alpha}}
|D\Phi(y)| dy.$$
Observe that
\begin{eqnarray}\label{3.2}
&|D\Phi(y)|=1+O(|y'|),\\[1mm]
 & |x-\Phi(y)|^2=|x-y+y-\Phi(y)|^2= |x_1-y_1|^2+|y'|^2+O(|y'|^4) \label{3.2-1}
\end{eqnarray}
and
\begin{eqnarray}\label{3.3}
\Big|\int_{ Q_{ \eta}}\frac{O(|y'|^2) }{|x-\Phi(y)|^{N+2\alpha}} |D\Phi(y)| dy\Big|\leq c_8,
\end{eqnarray}
where we use the fact that $|y-\Phi(y)|=O(|y'|^2)$.
Then we have that
 \begin{eqnarray*}
\tilde \BBE(x_1)&:=& {\rm p.v.} \int_{ Q_{ \eta}}\frac{x_1^{\tau}-y_1^{\tau} }{\left(|x_1-y_1|^2+|y'|^2+O(|y'|^4)\right)^{\frac{N+2\alpha}2}}(1+O(|y'|)) dy
\\&=& {\rm p.v.} \int_0^\eta(x_1^{\tau}-y_1^{\tau}) \int_{B^{N-1}_\eta(0)} \frac{ (1+O(|y'|)) }{\left(|x_1-y_1|^2+|y'|^2 )\right)^{\frac{N+2\alpha}2}} dy'dy_1
\\&=&{\rm p.v.} \int_0^\eta\frac{x_1^{\tau}-y_1^{\tau}}{|x_1-y_1|^{1+2\alpha}}   \int_{B^{N-1}_{\frac{\eta}{|x_1-y_1|}(0)}}
 \frac{ (1+O(|x_1-y_1| |z'|)) }{\left(1+|z'|^2 )\right)^{\frac{N+2\alpha}2}} dz' dy_1
 \\&=&{\rm p.v.} \int_0^\eta \frac{x_1^{\tau}-y_1^{\tau}}{|x_1-y_1|^{1+2\alpha}} F_\alpha(|x_1-y_1|)dy_1,
\end{eqnarray*}
 where
 \begin{eqnarray}
\nonumber && F_\alpha(|x_1-y_1|)
\\&=&\int_{B^{N-1}_{\frac{\eta}{|x_1-y_1|}(0)}}
 \frac{ (1+O(|x_1-y_1| |z'|)) }{\left(1+|z'|^2 )\right)^{\frac{N+2\alpha}2}} dz'
 \\\nonumber  &=&\int_{B^{N-1}_{\frac{\eta}{|x_1-y_1|}(0)}} \frac1{\left(1+|z'|^2 )\right)^{\frac{N+2\alpha}2}}dz'+O(|x_1-y_1|)
 \int_{B^{N-1}_{\frac{\eta}{|x_1-y_1|}(0)}} \frac{|z'|}{\left(1+|z'|^2 )\right)^{\frac{N+2\alpha}2}}dz'
 \\\nonumber &=&   \int_{\R^{N-1}} \frac1{\left(1+|z'|^2 )\right)^{\frac{N+2\alpha}2}}dz' +O(|x_1-y_1|^{1+2\alpha}) +\int_{\R^{N-1}} \frac{|z'|}{\left(1+|z'|^2 )\right)^{\frac{N+2\alpha}2}}dz' O(|x_1-y_1|)
 \\&=&d_{\alpha} +O(|x_1-y_1|)\label{3---1},
\end{eqnarray}
with $d_{\alpha}=\int_{\R^{N-1}} \frac1{\left(1+|z'|^2 )\right)^{\frac{N+2\alpha}2}}dz'>0$ and $B^{N-1}_r(0)$ being the ball
with center at the origin and radius $r$ in $\R^{N-1}$. Thus we have that
 \begin{eqnarray}
\nonumber \BBE(x_1)&=& d_{\alpha}\, {\rm p.v.} \int_0^\eta
 \frac{x_1^{\tau}-y_1^{\tau}}{|x_1-y_1|^{1+2\alpha}}  dy_1+O\Big({\rm p.v.} \int_0^\eta
 \frac{x_1^{\tau}-y_1^{\tau}}{|x_1-y_1|^{2\alpha}}  dy_1\Big)+O(1)
\\&=& \nonumber d_{\alpha} \left(x_1^{\tau-2\alpha} {\rm p.v.} \int_0^{\frac\eta{x_1}}
 \frac{1-t^{\tau}}{|1-t|^{1+2\alpha}}  dt+O(x_1^{\tau-2\alpha+1} {\rm p.v.} \int_0^{\frac\eta{x_1}}
 \frac{1-t^{\tau}}{|1-t|^{2\alpha}}  dt)\right)+O(1)
 \\&=& d_{\alpha} \Big[x_1^{\tau-2\alpha} \left(c_\alpha(\tau)-(\frac{\eta}{x_1})^{\tau-2\alpha}\right)+
 O(x_1^{\tau-2\alpha+1}\left(c_\alpha(\tau)-(\frac{\eta}{x_1})^{\tau-2\alpha+1}\right)\Big]+O(1)\nonumber
 \\&=&  d_{\alpha} c_\alpha(\tau)x_1^{\tau-2\alpha}+O(x_1^{\tau-2\alpha+1})+O(1).\label{E2(x1)+}
\end{eqnarray}

Therefore,  we conclude from   (\ref{E1(x1)+}) and  (\ref{E2(x1)+})
   that
\begin{eqnarray}
(-\Delta)^\alpha_\Omega V_\tau(x)= d_{\alpha} c_\alpha(\tau)x_1^{\tau-2\alpha} +O(x_1^\tau)+O(x_1^{\tau-2\alpha+1})+O(1)\quad{\rm as}\ \ {x_1\to0^+}.
\end{eqnarray}
Combining with Lemma \ref{lm 3.1},    Proposition \ref{pr 1.1} follows and the proof is complete.$\hfill\Box$\medskip

Note that $I_\alpha=\{0\}$ when $\alpha=\frac12$ and in this case, we need do the estimate for the following function
 \begin{equation}\label{1.1-log}
V_*(x)=\left\{ \arraycolsep=1pt
\begin{array}{lll}
 -\ln \rho(x) \ \ \ &\text{ for } \, x\in A_\delta, \\[2mm]
 l(x)\ \ \  &\text{ for } \, x\in \Omega\setminus A_\delta,
 \end{array}
\right.
\end{equation}
where $\delta\in(0,\frac14)$ and the function $l$  is  positive  such that
$V_*$ is $C^2$ in $\Omega$.

\begin{remark}\label{rm 3.1}
For $\tau\in(-1,2\alpha)\setminus\{2\alpha-1,0\}$, then
$$\lim_{\rho(x)\to0^+} \rho^{2\alpha-\tau}(x)(-\Delta)^\alpha_\Omega V_\tau(x)=d_{\alpha} c_\alpha(\tau).$$

\end{remark}

 \begin{proposition}\label{pr 1.1-ln}
 Assume that    $\Omega$ is a bounded $C^2$ domain in $\BBR^N$ with $N\geq 1$  and $V_*$ is given in (\ref{1.1-log}).

 Then there exists $c_9>0$  such that
  $$
|(-\Delta)^{1/2}_\Omega V_*(x)|\leq c_9 V_*(x), \ \ \ \ \forall\,x\in \Omega.$$
 \end{proposition}
\noindent{\bf Proof.} We omit the proof for the case $N=1$.  When $N\geq 2$ and under the same setting in the proof of Proposition \ref{pr 1.1}, we have
\begin{eqnarray}
(-\Delta)^{1/2}_\Omega V_*(x)=c_{N,\alpha}{\rm p.v.} \int_{\Phi(Q_{ \eta})}\frac{-\ln x_1+\ln \rho(z) }{|x-z|^{N+1}}dz+c_{N,\alpha} \int_{\Omega\setminus
\Phi(Q_{ \eta})}\frac{V_*(x)-V_*(z) }{|x-z|^{N+1}}dz \label{cotadelta-ln}
\end{eqnarray}
and we see that
\begin{equation}\label{E1(x1)+-ln}
|\int_{\Omega\setminus\Phi(Q_{ \eta})}\frac{V_*(x)-V_*(z) }{|x-z|^{N+2\alpha}}dz|
\le c_{10} V_*(x),
\end{equation}
where the constant $c_{10}$ is independent of $x$. Thus we only need to study the asymptotic behavior of the first integral, from now on, we denote by
$$\BBE(x_1)={\rm p.v.} \int_{\Phi(Q_{ \eta})}\frac{-\ln x_1+\ln \rho(z)  }{|x-z|^{N+1}}dz.$$
Let $z=\Phi(y)$ and $0<x_1<<1$, it follows from (\ref{3.2}), (\ref{3.2-1}) and (\ref{3.3}) that
\begin{eqnarray*}
\BBE(x_1)&=&{\rm p.v.} \int_{ Q_{ \eta}}\frac{-\ln x_1+\ln y_1 }{|x-\Phi (y)|^{N+1}}
|D\Phi(y)| dy+\int_{ Q_{ \eta}}\frac{O(|y'|^2) }{|x-\Phi (y)|^{N+1}}
|D\Phi(y)| dy,
\\&=& {\rm p.v.} \int_{ Q_{ \eta}}\frac{-\ln x_1+\ln y_1 }{\left(|x_1-y_1|^2+|y'|^2+O(|y'|^4)\right)^{\frac{N+1}2}}(1+O(|y'|)) dy
\\&=& {\rm p.v.} \int_0^\eta(-\ln x_1+\ln y_1) \int_{B^{N-1}_\eta(0)} \frac{ (1+O(|y'|)) }{\left(|x_1-y_1|^2+|y'|^2 )\right)^{\frac{N+1}2}} dy'dy_1
\\&=&{\rm p.v.} \int_0^\eta\frac{-\ln x_1+\ln y_1}{|x_1-y_1|^{1+1}}   \int_{B^{N-1}_{\frac{\eta}{|x_1-y_1|}(0)}}
 \frac{ (1+O(|x_1-y_1| |z'|)) }{\left(1+|z'|^2 )\right)^{\frac{N+1}2}} dz' dy_1
 \\&=&{\rm p.v.} \int_0^\eta \frac{-\ln x_1+\ln y_1}{|x_1-y_1|^2} F_{\frac12}(|x_1-y_1|)dy_1,
\end{eqnarray*}
 where  $F_{\frac12}(|x_1-y_1|)= d_{\frac12} +O(|x_1-y_1|)$, $d_{\frac12}=\int_{\R^{N-1}} \frac1{\left(1+|z'|^2 )\right)^{\frac{N+1}2}}dz'>0$.

Thus we have that
 \begin{eqnarray}
\nonumber \BBE(x_1)&=&d_{\frac12}\, {\rm p.v.} \int_0^\eta
 \frac{-\ln x_1+\ln y_1}{|x_1-y_1|^2}  dy_1+O({\rm p.v.} \int_0^\eta
 \frac{-\ln x_1+\ln y_1}{|x_1-y_1| }  dy_1)+O(1)
\\&=& \nonumber d_{\frac12} \left[x_1^{-1} {\rm p.v.} \int_0^{\frac\eta{x_1}}
 \frac{\ln t}{|1-t|^2}  dt+O\left(  \int_0^{\frac\eta{x_1}}
 \frac{\ln t }{|1-t| }  dt\right)\right]+O(1)
 \\&=& d_{\frac12} \Big[x_1^{-1} \left( -(\frac{\eta}{x_1})^{-1}\right)+
 O\left(\ln\ln \frac{\eta}{x_1}\right)+O(1)\Big] \nonumber
 \\&=&  O(\ln\frac1{x_1})+O(1),\label{E2(x1)+ln}
\end{eqnarray}
where we use (\ref{2.11}) and ${\rm p.v.} \int_0^{\frac\eta{x_1}}
 \frac{\ln t}{|1-t|^2}  dt= -\int_{\frac\eta{x_1}}^{+\infty}
 \frac{\ln t}{|1-t|^2}  dt$,
\begin{eqnarray*}
 \int_{\frac\eta{x_1}}^{+\infty} \frac{\ln t}{|1-t|^2}  dt &\leq &  2\int_{\frac\eta{x_1}}^{+\infty} \frac{\ln t}{t^2}  dt\\&=& -\frac{\ln t}{t}\Big|_{t=\frac{\eta}{x_1}} ^{+\infty}+\frac{1}{t}\Big|_{t=\frac{\eta}{x_1}} ^{+\infty}
\leq  -c_\eta x_1\ln x_1+c_\eta x_1.
\end{eqnarray*}

  Then we conclude from   (\ref{E1(x1)+-ln}) and  (\ref{E2(x1)+ln})
   that
\begin{eqnarray}
|(-\Delta)^{ 1/2}_\Omega V_*(x)|\leq c_{11}V_*(x)\quad\text{  for all}\  x\in\Omega,
\end{eqnarray}
which ends the proof.$\hfill\Box$

\mysection{Non-existence of positive solution }

\begin{proposition}\label{pr 3.1}
For $\alpha\in(0,\frac12]$, 
there is no   boundary blowing up regional fractional super-Harmonic function, i.e. 
 (\ref{eq 1.3}) has no super solution.
\end{proposition}
\noindent{\bf Proof. }  By contradiction, let (\ref{eq 1.3}) have a super solution $u$, then
its minimal is achieved inside, thanks to the boundary blowing up. Let $x_0$ be the point such that 
$$u(x_0)=\min_{x\in \Omega} u(x).$$ 
then in the viscosity sense,
\begin{equation}\label{y 2.1}
 (-\Delta)^{\alpha}_\Omega w(x_0)< 0,
\end{equation}
since $w$ blows up at the boundary, which contradicts the assumption that 
$w$ is a viscosity super solution of (\ref{eq 1.3}).\hfill$\Box$\smallskip

\begin{corollary}\label{cr 3.1}
For $\alpha\in(0,\frac12]$, 
problem
 \begin{equation}\label{eq 1.3-sub}
\arraycolsep=1pt\left\{
\begin{array}{lll}
 \displaystyle  (-\Delta)^\alpha_\Omega   u=0\quad & {\rm in}\ \,  \Omega,\\[1.5mm]
\phantom{ (-\Delta)^\alpha  }
 \displaystyle   u=-\infty\quad & {\rm on}\   \partial\Omega
 \end{array}\right.
\end{equation}
has no sub solution.
\end{corollary}

\begin{lemma}\label{lm 4.1}
Let
$$w_1(x)=V_*(x)-V_{\tau_0}(x)\quad{\rm for}\ \, x\in\Omega,$$
where $\tau_0=\frac14$.

Then there exist $\delta_2>0$ and $c_{12}>1$ such that
$$\frac1{c_{12}}\rho^{\tau_0-1}(x)\leq (-\Delta)^{1/2}_\Omega w_1(x)\leq c_{12}\rho^{\tau_0-1}(x),\quad\forall\, x\in A_{\delta_2}.$$
\end{lemma}
\noindent{\bf Proof. } From Proposition \ref{pr 1.1},  for $\delta_1>0$ and $c_5>1$,
$$\frac1{c_5}\rho^{-3/4}(x)\leq -(-\Delta)^{1/2}_\Omega V_{\tau_0}(x)\leq c_5\rho^{-3/4}(x),\quad\forall\, x\in A_{\delta_1}$$
and  from Proposition \ref{pr 1.1-ln},
$$|(-\Delta)^{1/2}_\Omega V_*(x)|\leq c_9V_*(x),\quad\forall\, x\in A_{\delta_1}.$$
Then there exist $\delta_2>0$ and $c_{12}>1$ such that
 $$\frac1{c_{12}}\rho^{\tau_0-1}(x)\leq (-\Delta)^{1/2}_\Omega w_1(x)\leq c_{12}\rho^{\tau_0-1}(x),\quad\forall\, x\in A_{\delta_2}.$$
The proof ends\hfill$\Box$ \medskip

    \noindent{\bf  Proof of Theorem \ref{teo 1}.}
    {\it Case 1:  $\alpha\in (0,\frac12)$. }  Take $\tau_0=\alpha-\frac12\in (2\alpha-1, 0)$ and from Proposition \ref{pr 1.1}, there exist $\delta_1>0$ and $c_{12}>0$ such that
   \begin{equation}\label{y 2.1-0}
   (-\Delta)^\alpha_\Omega V_{\tau_0}(x)\geq c_{12}\rho^{\tau_0-2\alpha},\quad x\in A_{\delta_1}. 
   \end{equation}
   Since $V_{\tau_0}$ is $C^2$ in $\Omega$, then there exists $T_0>0$ such that
     \begin{equation}\label{y 2.1-1}
     |(-\Delta)^\alpha_\Omega V_{\tau_0}(x)|\leq T_0\quad {\rm for}\ \, x\in \Omega\setminus A_{\delta_1}.   
     \end{equation}

 By contradiction, we assume that problem (\ref{eq 1.1})
    has a  solution $u_1$ bounded from below.

     Let
    $$w=V_{\tau_0}+t u_1 \quad{\rm in}\ \ \Omega,$$
    where $t>0$ will be choosing later. 

In the viscosity sense, we have that
\begin{equation}\label{y 2.1-2}
 (-\Delta)^{\alpha}_\Omega w= (-\Delta)^{\alpha}_\Omega V_{\tau_0}+t   \geq 0\quad{\rm in}\ \, \Omega
\end{equation}
if $t\geq T_0$. 

Therefore,  $w$ is a viscosity super-solution of $(-\Delta)^{\alpha}_\Omega u=0$ and  blows up at boundary of $\Omega$, since $u_1$ is bounded from below.  A contradiction 
comes from Proposition \ref{pr 3.1}. 

 By contradiction, we assume that problem (\ref{eq 1.1})
    has a  solution $u_1$ bounded from above.   Let
    $$w=-V_{\tau_0}+t u_1 \quad{\rm in}\ \ \Omega $$
 for $t<0$ suitable, then $w$ will be a solution of (\ref{eq 1.3-sub}), which is impossible. 

\smallskip

 {\it Case 2:  $\alpha=\frac12$.} From Lemma \ref{lm 4.1}, we have that
     \begin{equation}\label{y 2.2-0}
     (-\Delta)^\alpha_\Omega w_1(x)\geq c_{12}\rho^{\tau_0-2\alpha},\quad x\in A_{\delta_1}.
     \end{equation}
   Since $w_1$ is $C^2$ in $\Omega$, then there exists $T_1>0$ such that
       \begin{equation}\label{y 2.2-1}
       |(-\Delta)^\alpha_\Omega V_{\tau_0}(x)|\leq T_1\quad {\rm for}\ \, x\in \Omega\setminus A_{\delta_1}.     
       \end{equation}
   The rest proof is the same as {\it Case 1}. \hfill$\Box$\medskip

   \noindent{\bf Proof of Theorem \ref{teo 2}. }  Assume that $u_p\geq 0$ is a nontrivial classical solution of  (\ref{eq 1.2}),
   if there exists $x_0\in\Omega$ such that $u_p(x_0)=0$, then it is the minimum point in $\Omega$ and the definition of fractional laplacian implies
that
   $$0>(-\Delta)^\alpha_\Omega u_p(x_0)=u_p^p(x_0)= 0,$$
 which is impossible. 
 Thus, we have that $u_p>0$ in $\Omega$. 
 
 For $\delta_1>0$ in (\ref{y 2.2-0}),  there exists $\varepsilon_0>0$ such that 
 $$u_p \geq \varepsilon_0\quad {\rm in}\ \, \Omega\setminus A_{\delta_1}.$$
 
  Let
    $$w=V_{\tau_0}+t u_p\quad{\rm in}\ \ \Omega,$$
    where $t>0$ will be choosing later. 

In the viscosity sense, we have that
$$
 (-\Delta)^{\alpha}_\Omega w= (-\Delta)^{\alpha}_\Omega V_{\tau_0}+t\varepsilon_0^{p}   \geq 0\quad{\rm in}\ \, \Omega\setminus A_{\delta_1}.
$$
If $t\geq T_0\varepsilon_0^{-p}$, then 
  $w$ is a viscosity super-solution of $(-\Delta)^{\alpha}_\Omega u=0$ and  blows up at boundary of $\Omega$.   A contradiction 
comes from Proposition \ref{pr 3.1}.  
 \hfill$\Box$

\section*{ Appendix: properties of killing measure }

The killing measure $\kappa_\alpha$ is the connection between global fractional Laplcian and regional fractional Laplcian. For a regular function $w$ such that $w=0$ in $\BBR^N\setminus\bar\Omega$, we observe that
  $$
   (-\Delta)^\alpha_\Omega w(x)= (-\Delta)^\alpha w(x)-w(x)\kappa_\alpha(x),\quad\forall x\in\Omega,
 $$
where
$$
 \kappa_\alpha(x)=c_{N,\alpha}\int_{\BBR^N\setminus \Omega} \frac1{|x-z|^{N+2\alpha}}dz.
$$
\smallskip

 \noindent{\bf Proposition A.1}  {\it
Let $N\geq 2$, $\Omega$be a $C^2$ domain  and $\rho(x)=dist(x,\partial\Omega)$, then   $\kappa_\alpha\in C^{0,1}_{\rm loc}(\Omega)$ and
\begin{equation}\label{2.6}
 \lim_{\rho(x)\to0^+}\kappa_\alpha(x)\rho(x)^{2\alpha}=d_\alpha c_{N,\alpha},
\end{equation}
where
$$
 d_\alpha=\int_{\R^{N-1}}  \frac{1}{ (1+|z'|^2)^{\frac{N+2\alpha}{2}}} dz'.
$$
}

\noindent{\bf Proof.} {\it Lipchitz continuity of $\kappa_\alpha$. } For
 $x_1,x_2\in \Omega $ and  any $z\in \BBR^N\setminus \Omega$, we have that
 $$|z-x_1|\geq \rho(x_1)+\rho(z), \qquad |z-x_2|\geq \rho(x_2)+\rho(z)$$ and
$$||z-x_1|^{N+2\alpha}-|z-x_2|^{N+2\alpha}|\leq
c_{13}(|x_1-x_2|(|z-x_1|^{N+2\alpha-1}+|z-x_2|^{N+2\alpha-1}),
$$
for some $c_2>0$ independent of $x_1$ and $x_2$. Then
\begin{eqnarray*}
&&|\kappa_\alpha(x_1)-\kappa_\alpha(x_2)| \leq c_{N,\alpha}\int_{\R^N\setminus \Omega}  \frac{||z-x_2|^{N+2\alpha}-|z-x_1|^{N+2\alpha}|}{|z-x_1|^{N+2\alpha}|z-x_2|^{N+2\alpha}}dz
\\&&\leq c_{N,\alpha} c_{13}|x_1-x_2|\left[\int_{\R^N\setminus \Omega}  \frac{dz}{|z-x_1||z-x_2|^{N+2\alpha}}+\int_{\BBR^N\setminus \Omega}  \frac{dz}{|z-x_1|^{N+2\alpha}|z-x_2|}\right].
\end{eqnarray*}
By direct computation, we have that
\begin{eqnarray*}
 \int_{\BBR^N\setminus \Omega}  \frac{1}{|z-x_1||z-x_2|^{N+2\alpha}}dz
  &\le& \int_{\BBR^N\setminus{B_{\rho(x_1)}(x_1)}} \frac{1}{|z-x_1|^{N+2\alpha+1}}dz
 \\&&+\int_{\BBR^N\setminus{B_{\rho(x_2)}(x_2)}} \frac{1}{|z-x_2|^{N+2\alpha+1}}dz \\
 &\le& c_{14}[\rho(x_1)^{-1-2\alpha}+\rho(x_2)^{-1-2\alpha}]
\end{eqnarray*}
and similar to obtain that
$$\int_{\BBR^N\setminus \Omega} \frac{1}{|z-x_1|^{N+2\alpha}|z-x_2|}dz\le  c_{15}[\rho(x_1)^{-1-2\alpha}+\rho(x_2)^{-1-2\alpha}],$$
where $c_{15}>0$ is independent of $x_1,\, x_2$.
Then
$$|\kappa_\alpha(x_1)-\kappa_\alpha(x_2)|\le c_{14} c_{15}[\rho(x_1)^{-1-2\alpha}+\rho(x_2)^{-1-2\alpha}]|x_1-x_2|,$$
that is,
 $\kappa_\alpha$ is $C^{0,1}$ locally in $\Omega$.\smallskip

 {\it Proof of (\ref{2.6}).}    By compactness we prove that the corresponding inequality holds in a neighborhood of any point $\bar x\in\partial\Omega$ and  without loss of generality we may assume that $\bar
x=0$ and we further assume that  $-e_1$ is the outer normal
vector of $\Omega$ at $0$, where $e_1=(1,0,\cdots)\in\BBR^N$.

 For  a given $0<\eta\le \delta$, we define
$$
Q_\eta^-=\{y=(y_1,y')\in\R\times\R^{N-1},\, -\eta< y_1\leq 0,\, |y'|<\eta\}$$
and
$$
\tilde Q_\eta=\{y=(y_1,y')\in\R\times\R^{N-1},\, |y_1|<\eta,\, |y'|<\eta\}.$$

For  $C^2$ domain $\Omega$,  we recall  $C^2$ diffeomorphism mapping $\Phi:\BBR^N\to \BBR^N$ that
$$\Phi(z)=(z', z_1+\varphi(z'))\quad{\rm for}\ \, z=(z_1,\, z')\in   Q_\eta.$$
 More properties of $\Phi$ is given in the proof of {\it Proposition \ref{pr 1.1}}.

Let $x_1\in(0,\eta/4)$,
$x=(x_1,0)$   be  a generic point in $A_{\eta/4}$.
 By the definition of $\kappa_\alpha$, we have that
 \begin{eqnarray*}
\frac1{c_{N,\alpha}}\kappa_\alpha(x)= \int_{\Phi(Q_\eta^-)} \frac1{|x-y|^{N+2\alpha}}dy  +\int_{ \Omega^c\setminus \Phi(Q_\eta^-)} \frac1{|x-y|^{N+2\alpha}}dy ,
 \end{eqnarray*}
 where $\Omega^c=\BBR^N\setminus \Omega$.

 On the one hand, we have that
 $  |x-z| \geq \frac14(\eta +|z|) \text{ for any } z\in \Omega^c\setminus \Phi(Q_\eta^-)$ and
 then
 $$\int_{ \Omega^c\setminus \Phi(Q_\eta^-)} \frac1{|x-z|^{N+2\alpha}}dz\leq \int_{\BBR^N} \frac{4^{N+2\alpha}}{(\eta+|z|)^{N+2\alpha}}dy \leq  c_{15},$$
 where $c_{15}>0$ depends on $\eta$.

  On the other hand, by change variable $z=\Phi(y)$,   we have that
$$ \int_{\Phi(Q_\eta^-)} \frac1{|x-z|^{N+2\alpha}}dz = \int_{ Q_{ \eta}^-}\frac{ |D\Phi(y)|  }{|x-\Phi(y)|^{N+2\alpha}}
dy .$$
Observe that
\begin{eqnarray*}
&|D\Phi(z)|=1+O(|z'|),\\[1mm]
 & |x-\Phi(y)|^2=|x-z+z-\Phi(y)|^2= |x_1-y_1|^2+|y'|^2+O(|y'|^4),
\end{eqnarray*}
then
\begin{eqnarray*}
\int_{\Phi(Q_\eta^-)} \frac1{|x-z|^{N+2\alpha}}dx& =& \int_{ Q_{ \eta}^-}\frac{ 1+O(|y'|)}{(|x_1-y_1|^2+|y'|^2)^{\frac{N+2\alpha}2}}
dz
\\&=& \int^0_{-\eta} \frac1{|x_1-y_1|^{1+2\alpha}}  \int_{B^{N-1}_{\frac{\eta}{|x_1-y_1|}(0)}} \frac{1+O(|x_1-y_1||\zeta'|)}{\left(1+|\zeta'|^2 \right)^{\frac{N+2\alpha}2}}d\zeta' dz_1
\\&=&\frac1{2\alpha} \int_{\R^{N-1}}  \frac{1}{ (1+|\zeta'|^2)^{\frac{N+2\alpha}{2}}} d\zeta' \Big(x_1^{-2\alpha}-(x_1+\eta_1)^{-2\alpha}\Big)
\\&&+  O\Big(\int_{\R^{N-1}}  \frac{|\zeta'|}{ (1+|\zeta'|^2)^{\frac{N+2\alpha}{2}}} d\zeta' \varphi(x_1)\Big),
 \end{eqnarray*}
where $B^{N-1}_{r}(0)$ is the ball center at origin with radius $r$ in $\BBR^{N-1}$,
$\varphi_\alpha(t)=t^{1-2\alpha}$ if $2\alpha\not=1$ and $\varphi_{\frac12}(t)=-\log t$. The proof ends.\hfill $\Box$ \medskip\smallskip

 \noindent{\bf Remark A.1.} {\it
When the domain has the property of  uniformly interior tangent ball, i.e.
then it has the rough estimate:
$$\frac1{c_{16}}\rho(x)^{-2\alpha}\leq \kappa_\alpha(x)\leq c_{16}\rho(x)^{-2\alpha},\quad\forall\, x\in\Omega,$$
where $c_{16}>1$.
 }


\end{document}